\documentclass[12pt]{amsart}
\usepackage{amscd,amssymb}
\usepackage[graph,frame,poly,arc]{xy}  
\usepackage[plainpages,backref,urlcolor=blue]{hyperref}

\theoremstyle{plain}
\newtheorem{thm}[subsection]{Theorem}

\newtheorem{prop}[subsection]{Proposition}
\newtheorem{cor}[subsection]{Corollary}

\theoremstyle{definition}
\newtheorem{rk}[subsection]{Remark}
\newtheorem{definition}[subsection]{Definition}
\newtheorem{ex}[subsection]{Example}
\newtheorem{conj}[subsection]{Conjecture}

\numberwithin{equation}{section}
\newcommand{\OO}{{\mathcal O}}

\newcommand{\A}{{\mathcal A}}

\newcommand{\C}{\mathbb{C}}
\newcommand{\PP}{\mathbb{P}}

\DeclareMathOperator{\im}{im}

\DeclareMathOperator{\reg}{reg}

\begin{document}

\title [Curves with maximal global Tjurina numbers]
{Jacobian syzygies and plane curves with maximal global Tjurina numbers}

\author[Alexandru Dimca]{Alexandru Dimca$^{1}$}
\address{Universit\'e C\^ ote d'Azur, CNRS, LJAD and INRIA, France and Simion Stoilow Institute of Mathematics,
P.O. Box 1-764, RO-014700 Bucharest, Romania}
\email{dimca@unice.fr}

\author[Gabriel Sticlaru]{Gabriel Sticlaru}
\address{Faculty of Mathematics and Informatics,
Ovidius University
Bd. Mamaia 124, 900527 Constanta,
Romania}
\email{gabrielsticlaru@yahoo.com }

\thanks{$^1$ This work has been partially supported by the French government, through the $\rm UCA^{\rm JEDI}$ Investments in the Future project managed by the National Research Agency (ANR) with the reference number ANR-15-IDEX-01 and by the Romanian Ministry of Research and Innovation, CNCS - UEFISCDI, grant PN-III-P4-ID-PCE-2016-0030, within PNCDI III}



\begin{abstract}  First we give a sharp upper bound for the cardinal $m$ of a minimal set of generators for the module of Jacobian syzygies of a complex projective reduced plane curve $C$. Next we discuss the sharpness of an upper bound,
given by A. du Plessis and C.T.C. Wall,
for the global Tjurina number of such a  curve $C$, in terms of its degree $d$ and of the minimal degree $r\leq d-1$ of a Jacobian syzygy. We give a  homological characterization of the curves whose global Tjurina number equals the du Plessis-Wall upper bound, which implies in particular that for such curves the upper bound for $m$ is also attained.
Finally we prove the existence of curves with maximal global Tjurina numbers for certain  pairs $(d,r)$. Moreover, we conjecture that such curves exist for any pair $(d,r)$, and that, in addition, they may be chosen to be line arrangements when $r\leq d-2$. This conjecture is proved for degrees $d \leq 11$.

\end{abstract}
 
\maketitle


\section{Introduction} 

Let $S=\C[x,y,z]$ be the polynomial ring in three variables $x,y,z$ with complex coefficients, and let $C:f=0$ be a reduced curve of degree $d$ in the complex projective plane $\PP^2$. 
We denote by $J_f$ the Jacobian ideal of $f$, i.e. the homogeneous ideal in $S$ spanned by the partial derivatives $f_x,f_y,f_z$ of $f$, and  by $M(f)=S/J_f$ the corresponding graded quotient ring, called the Jacobian (or Milnor) algebra of $f$.
Consider the graded $S$-module of Jacobian syzygies of $f$, namely
$$AR(f)=\{(a,b,c) \in S^3 \ : \ af_x+bf_y+cf_z=0\}.$$
We say that $C:f=0$ is an {\it $m$-syzygy curve} if the graded $S$-module $AR(f)$ is generated by $m$ homogeneous syzygies, say $\rho_1,\rho_2,...,\rho_m$, with $m$ minimal possible, of degrees $d_j=\deg \rho_j$ ordered such that $$1 \leq d_1\leq d_2 \leq ...\leq d_m.$$ 
In fact, the case $d_1=0$ occurs if and only if $C$ is a union of lines through a point, a situation which is not considered in the sequel.
 We call these degrees the {\it exponents} of the curve $C$ and the syzygies $\rho_1,...,\rho_m$ a {\it minimal set of generators } for the module  $AR(f)$. Note that $d_1=mdr(f)$ is the minimal degree of a non trivial Jacobian syzygy in $AR(f)$.  The minimal possible value for $m$ is two, since the $S$-module $AR(f)$ has rank two for any $f$. The curve $C$ is called {\it free} when $m=2$, since then  $AR(f)$ is a free module, see for such curves \cite{B+,Dmax,DStFD,DStMos,Sim1,Sim2,ST,To}. 
 We prove in Proposition \ref{propNS}, which is our first result, that
for a reduced, degree $d$ curve $C$, one has 
\begin{equation}
\label{m1}
m \leq d_1+d_2-d+3\leq d+1.
\end{equation}
Note that for the case $C$ a line arrangement,  the slightly stronger inequality
 $m \leq d-1$ was known,
see \cite[Corollary 1.3]{DStSat}. 
Moreover, all these inequalities, involving the integers $m$, $d_1$, $d_2$ and $d$, are sharp, see Remark \ref{rkS1}.

 Recall that the {\it global Tjurina number} $\tau(C)$ of the plane curve $C:f=0$ can be defined as either the degree of the Jacobian ideal $J_f$, or as the sum of the Tjurina numbers of all the singularities of the curve $C$. With this notation, it was shown by A. du Plessis and C.T.C. Wall that one has the following result, see \cite[Theorem 3.2]{dPW}, and also \cite[Theorem 20]{E} for a new approach.

\begin{thm}
\label{thm1}
Let $C:f=0$ be a reduced plane curve of degree $d$ and let $r=d_1$  be the minimal degree of a non-zero syzygy in  $AR(f)$.
Then the following hold.
\begin{enumerate}

\item If $r <d/2$, then $\tau(C) \leq \tau(d,r)_{max}= (d-1)(d-r-1)+r^2$ and the equality holds if and only if the curve $C$ is free.

\item If $d/2 \leq r \leq d-1$, then
$\tau(C) \leq \tau(d,r)_{max}$,
where, in this case, we set
$$\tau(d,r)_{max}=(d-1)(d-r-1)+r^2-{ 2r-d+2 \choose 2}.$$

\end{enumerate}

\end{thm}
In this note we investigate for which curves one has equality in the above result. To have a name, we call such curves {\it maximal Tjurina curves of type} $(d,r)$.
Note that for any pair
$(d,r)$, with $1 \leq r <d/2$, the existence of maximal Tjurina curves of type $(d,r)$, i.e. of free curves with these invariants, follows from \cite{DStExpo}. 
The characterization and the existence of maximal Tjurina curves of type $(d,r)$, with $d/2 \leq r \leq d-1$,
 is our main concern in this note, and hence {\it we assume from now on that} $d/2\leq r$.  
 In the third section we derive a {\it homological characterization} of these maximal Tjurina curves, see Theorem \ref{thm2}. If we set $m=2r-d+3$, this result says  that a maximal Tjurina curve of type $(d,r)$ is exactly an $m$-syzygy curve, with exponents
 $$d_1=d_2= \cdots =d_m=r.$$
 In particular, a maximal Tjurina curve has the largest number $m$ of generators for the module $AR(f)$ allowed by the first inequality in \eqref{m1}. Another characterization of maximal Tjurina curves is given in
 \cite[Theorem 4.1]{DStFit}.
 In the last section we describe  some existence results for such curves.
 For $m=3$, the minimal possible value, a maximal Tjurina curve of type $(d,r)$ is exactly a
{\it nearly free} curve of degree $d=2r$, with exponents $d_1=d_2=d_3=r$, see subsection \ref{s1},
and hence the existence of maximal Tjurina curves for any type $(d,r)=(2r,r)$ follows again from \cite{DStExpo}. 

For $m=4$, a maximal Tjurina curve of type $(d,r)$ is a 4-syzygy curve with $d=2r-1$ and exponents $d_1=d_2=d_3=d_4=r$. These curves are related to {\it nearly cuspidal rational curves}, i.e. to rational curves having only unibranch singularities, except from one singularity which has 2 branches, see subsection \ref{s2}. It is conjectured that a nearly cuspidal rational curve $C$ satisfies the inequalities $m\leq 4$ and
$$\tau(d,r)_{max}-2 \leq \tau(C) \leq \tau(d,r)_{max},$$
see \cite{DSt3syz}, and moreover $\tau(C) = \tau(d,r)_{max}$ when $m=4$.
Proposition \ref{prop0} below describe a sequence of nearly cuspidal rational curves $C_d$ which are  maximal Tjurina curves of type $(2r-1,r)$ for any $r\geq 3$.

At the other extreme, for $m=d+1$, the maximal possible value for $m$, we notice that in this case $r=d-1$ and an example of maximal Tjurina curve of type $(d,d-1)$ is given by any maximal nodal rational curve of degree $d$, see Proposition \ref{prop1}.

If we go one step back, for $m=d-1$, we have $r=d-2$ and an example of maximal Tjurina curve of type $(d,d-2)$ is given by any generic  arrangement of $d$ lines in $\PP^2$, see
Proposition \ref{prop0.1}.
If we go back one more step, namely for $m=d-3$, and hence $r=d-3$, we describe a sequence of line arrangements $C_d$, having only double and triple intersection points, which are very likely maximal Tjurina curves of type 
$(d,d-3)$ for any $d \geq 6$, and we check this claim for $d \in [6,20]$ using SINGULAR, see \cite{Sing}. Similarly, for $m=d-5$ and $r=d-4$
we describe a sequence of line arrangements $\A_d$, having only double, triple and 4-fold intersection points, which are very likely maximal Tjurina curves of type 
$(d,d-4)$ for any $d \geq 8$, and we check this claim for $d \in [8,20]$ using SINGULAR, see \cite{Sing}.


In view of all these examples, we offer the following.

\begin{conj}
\label{conj1}  For any integer $d \geq 3$ and for any integer $r$ such that $d/2\leq r \leq d-1$,
there are maximal Tjurina curves of type $(d,r)$. Moreover, for $d/2 \leq r\leq d-2$, there are maximal Tjurina line arrangements of type $(d,r)$.
\end{conj}
In other words, the du Plessis-Wall inequality in Theorem \ref{thm1} is sharp for any pair
$(d,r)$ as above. The fact that line arrangements seem to give examples of maximal Tjurina curves of type $(d,r)$ for any $r<d-1$ may encourage further study of the deep relation between the combinatorics of a line arrangement $\A:f=0$ and the integer $r=mdr(f)$, see the end of subsection \ref{s3} for a brief discussion of this point.
The existence results in the final section imply the following.
\begin{cor}
\label{corCONJ}  Conjecture \ref{conj1} holds for $d \leq 11$.
\end{cor}
 

We would like to thank Aldo Conca for useful discussions.

\section{On the number of generators of the Jacobian syzygies} 

Consider the general form of the minimal graded resolution for the graded $S$-module 
  $M(f)$, the Milnor algebra of a curve $C:f=0$, that is assumed to be not free, namely
\begin{equation}
\label{res2A}
0 \to \oplus_{i=1} ^{m-2}S(-e_i) \to \oplus_{j=1} ^mS(1-d-d_j)\to S^3(1-d)  \to S,
\end{equation}
with $m\geq 3$, $e_1\leq ...\leq e_{m-2}$ and $d_1\leq ...\leq d_m$.
Since the kernel of the map $S^3(1-d)  \to S$ is precisely $AR(f)(1-d)$, we see that $d_1\leq ...\leq d_m$ are the exponents of $C$ as in the Introduction.
It follows from \cite[Lemma 1.1]{HS} that one has
\begin{equation}
\label{res2B}
e_j=d+d_{j+2}-1+\epsilon_j,
\end{equation}
for $j=1,...,m-2$ and some integers $\epsilon_j\geq 1$. Using  \cite[Formula (13)]{HS}, it follows that one has
\begin{equation}
\label{res2C}
d_1+d_2=d-1+\sum_{i=1} ^{m-2}\epsilon_j.
\end{equation}
It is known that, for a reduced degree $d$ curve $C$, one has $$d_m \leq 2d-4,$$
 see \cite[Corollary 11]{CD}, as well as  \cite[Corollary 12]{E3} for the case of a quasi-complete intersection ideal replacing $J_f$ and also 
 \cite[Theorem 9.4]{Chardin} for an even more general case. For the case $C$ a line arrangement, one has the much stronger inequality
 $d_m \leq d-2,$
see \cite[Corollary 3.5]{Sch0}. The inequality $d_m \leq d-1$  holds more generally for curves $C$ having as irreducible components only rational curves, see \cite[Corollary 5.2]{DSt3syz}.
The first main result of this note is the following.
\begin{prop}
\label{propNS}
Let $C:f=0$ be an $m$-syzygy curve of degree $d\geq 3$,
with exponents $1\leq d_1 \leq d_2\leq \cdots  \leq d_m$. Then $m$, the cardinal of a minimal set of generators for the module  $AR(f)$,
 satisfies the inequalities
$$m \leq d_1+d_2-d+3 \leq d_1+2 \leq d+1.$$
\end{prop}
This result was obtained independently by Philippe Ellia,
see Corollary 5 (i) and Theorem 7 (ii) in \cite{E3}, with a different approach and in a more general setting: the Jacobian ideal $J_f$ is replaced by a quasi-complete intersection.
\proof
If $m=2$, then the curve $C$ is free, and there is nothing to prove,
since for such curves $d_1+d_2=d-1$. 
Assume that $m \geq 3$. The first claim follows from the equality \eqref{res2C}, which yields $m-2 \leq d_1+d_2-d+1$,
since all the numbers $\epsilon_j$ are strictly positive integers.
To get the other two inequalities, recall that  \cite[Theorem 2.4]{DSt3syz} implies that
$$d_1\leq d_2 \leq d-1.$$
 \endproof
\begin{rk}
\label{rkS1}
 (i) One has the equality $m=d+1$ in Proposition \ref{propNS} for some curves, in particular for any maximal nodal rational curve, see Example \ref{ex1}, Example \ref{ex2} and Proposition \ref{prop1} below.  Moreover, the equality
 $m =d_1+d_2-d+3$ holds for any maximal Tjurina curve of type $(d,r)$, with
 $r=mdr(f)\geq d/2$, as shown in Theorem \ref{thm2} below. In this case $d_1=d_2=r$.

\medskip

\noindent (ii) Recall that for a generic, i.e. nodal, line arrangement $C$ in $\PP^2$, one has $m=d-1$ and $d_m=d-2$, see
\cite[Corollary 3.5]{Sch0}.
The fact that for any line arrangement $C$, one has $m \leq d-1$, see \cite[Corollary 1.3]{DStSat}, can in fact be proven using the same
idea as in the proof of  Proposition \ref{propNS}. Indeed, if the the arrangement is not generic, it follows that it has a point of multiplicity $m \geq 3$. Then \cite[Theorem 1.2]{Dcurves} implies that 
$$d_1 \leq d-m \leq d-3,$$
and hence 
$d_1 +d_2 \leq 2d-4$, since $d_2 \leq d-1$ by  \cite[Theorem 2.4]{DSt3syz}.
This yields $m \leq d-1$ in this case.

\medskip

\noindent (iii) Note that, for a uninodal plane curve of degree $d$, the module $AR(f)$ has 4 minimal generators, of degrees $d-1, d-1, d-1, 2d-4$.  Such curves provide also examples for the equality $d_m=2d-4$.

\end{rk}

\section{A characterization of maximal Tjurina curves} 

We recall now the construction of the Bourbaki ideal $B(C,\rho_1)$ associated to a degree $d$ reduced curve $C:f=0$ and to a minimal degree non-zero syzygy $\rho_1 \in AR(f)$, see \cite{DStJump}.
For any choice of the syzygy $\rho_1=(a_1,b_1,c_1)$ with minimal degree $r=d_1$, we have a morphism of graded $S$-modules
\begin{equation} \label{B1}
S(-r)  \xrightarrow{u} AR(f), \  u(h)= h \cdot \rho_1.
\end{equation}
For any homogeneous syzygy $\rho=(a,b,c) \in AR(f)_m$, consider the determinant $\Delta(\rho)=\det M(\rho)$ of the $3 \times 3$ matrix $M(\rho)$ which has as first row $x,y,z$, as second row $a_1,b_1,c_1$ and as third row $a,b,c$. Then it turns out that $\Delta(\rho)$ is divisible by $f$, see \cite{Dmax}, and we define thus a new morphism of graded $S$-modules
\begin{equation} \label{B2}
 AR(f)  \xrightarrow{v}  S(r-d+1)   , \  v(\rho)= \Delta(\rho)/f,
\end{equation}
and a homogeneous ideal $B(C,\rho_1) \subset S$ such that $\im v=B(C,\rho_1)(r-d+1)$.
It is known that the ideal $B(C,\rho_1)$, when $C$ is not a free curve, defines a $0$-dimensional subscheme $Z(C,\rho_1)$ in $\PP^2$, which is locally a complete intersection, see \cite[Theorem 5.1]{DStJump}. Using this construction, we can prove the following characterization of maximal Tjurina curves, which is our second main result in this paper.
\begin{thm}
\label{thm2}

Let $C:f=0$ be a reduced plane curve of degree $d$, let $r=mdr(f)$  be the minimal degree of a non-zero syzygy in  $AR(f)$ and assume 
 $d/2 \leq r \leq d-1$. Then
$\tau(C) \leq \tau(d,r)_{max}$,
and, if  equality holds, then the minimal resolution of the graded $S$-module AR(f) has the form
$$0 \to S(-r-1)^{m-2} \to S(-r)^m \to AR(f) \to 0,$$
where $m=2r-d+3$.
In particular, the exponents of the curve $C$ are given by
$$d_1=d_2= \cdots = d_m=r.$$
Conversely, if $C:f=0$ is a reduced plane curve of degree $d$, which has exponents
$$d_1=d_2= \cdots = d_m=r,$$
with $m=2r-d+3$, then the curve $C:f=0$ is a maximal Tjurina curve of type $(d,r)$.

\end{thm}
\proof Since the quotient $S^3/AR(f)$ is torsion free, it follows that the ideal $I=B(C,\rho_1)$ is saturated, and hence $P=S/I$ is a Cohen-Macaulay module. This fact has two consequences. First the Hilbert function $H_P(k)= \dim P_k$ is increasing. By definition, all the generators of $I$ have degree at least $2r-d+1$, and hence we get
\begin{equation} \label{B3}
{ 2r-d+2 \choose 2}=\dim S_{2r-d}= \dim P_{2r-d} \leq H_P(k),
\end{equation}
for large $k$. On the other hand, on has 
\begin{equation} \label{B4}
H_P(k) =\deg Z(C,\rho_1)=(d-1)^2-r(d-r-1)-\tau(C),
\end{equation}
for large $k$, see \cite[Theorem 5.1]{DStJump}. The last two relations imply the du Plessis-Wall inequality. Moreover, we see that we have equality for the curve $C$ if and only if
\begin{equation} \label{B5}
{ 2r-d+2 \choose 2}= \dim P_{2r-d}=  H_P(k),
\end{equation}
for all $k \geq 2r-d$. Since $P$ is a Cohen-Macaulay module, it follow that
$$\reg (P)=2r-d,$$
where $\reg(P)$ denotes the Castelnuovo-Mumford regularity of the $S$-module $P$, see
\cite[Theorem 4.2]{Eis}.
The minimal resolution of $P$ has the form
$$ 0 \to \oplus_jS(-a_{2,j}) \to \oplus_jS(-a_{1,j}) \to S \to P \to 0,$$
where $a_{1,j} \geq 2r-d+1$ are the degrees of the generators for the ideal $I$. It follows that all these generators must have degree $a_{1,j}=2r-d+1$, since by definition
$$\reg(P)=\max_{i,j}(a_{i,j}-1).$$
In order to have $H_P(2r-d)=H_P(2r-d+1)$, we need exactly
$$m'= \dim S_{2r-d+1}- \dim S_{2r-d}=2r-d+2$$
generators for $I$. It follows that the above minimal resolution for $P$ yields the following minimal resolution
$$0 \to S(d-2-2r)^{m'-1} \to S(d-1-2r)^{m'} \to I \to 0,$$
for the ideal $I=B(C,\rho_1)$.
Using the exact sequence 
$$0 \to S(-r) \to AR(f) \to B(C,\rho_1)(r-d+1) \to 0,$$
it follows that $AR(f)$ is minimally generated by $m=m'+1$ generators, all of degree $r$,
the first one being $\rho_1$, and then $\rho_j$ for $j=2,...,m$ being chosen such that their images under $v$ generate the ideal $I$. Moreover, each of the $m'-1$ relations among the generators of $I$ will give rise to a relation, with linear coefficients, among the syzygies
$\rho_i$. It follows that the minimal resolution of the $S$-module $AR(f)$ is given by
$$0 \to S(-r-1)^{m-2} \to S(-r)^m \to AR(f) \to 0.$$
To prove the converse, it is enough to show that our hypothesis implies that the minimal resolution of the $S$-module $AR(f)$ has the form above.
Indeed, the minimal resolution of the $S$-module $AR(f)$ determines
both $r=mdr(f)$ and $\tau(C)$, e.g. using the exact sequence in Corollary \ref{cor2} below.
To show that the minimal resolution of the $S$-module $AR(f)$ has the form above, we use the formula \eqref{res2C}.
This  formula  implies that $\epsilon_j=1$ for any $j$, and hence all the second order syzygies of
$AR(f)$ have the same degree
$$e'_1=e'_2= \cdots =e'_{m-2}=r+1.$$
This implies that $e_j=e_j'+d-1=r+d$, for all $1 \leq j \leq m-2$, which
 ends the proof of Theorem \ref{thm2}.
\endproof

\begin{rk}
\label{rkUlrich}  If we denote by $E_C$ the rank two vector bundle on $\PP^2$ associated to the graded $S$-module $AR(f)$, then Theorem \ref{thm2} implies the existence of an exact sequence
$$0 \to \OO_{\PP^2}(-r-1)^{m-2} \to \OO_{\PP^2}(-r)^{m} \to E_C \to 0,$$
for any maximal Tjurina curve $C$.
\end{rk}

Recall the following definition, see \cite{Dmax, DStFD}.

\begin{definition}
\label{def}

\noindent (i) the {\it coincidence threshold} 
$$ct(f)=\max \{q:\dim M(f)_k=\dim M(f_s)_k \text{ for all } k \leq q\},$$
with $f_s$  a homogeneous polynomial in $S$ of the same degree $d$ as $f$ and such that $C_s:f_s=0$ is a smooth curve in $\PP^2$.

\noindent (ii) the {\it stability threshold} 
$st(f)=\min \{q~~:~~\dim M(f)_k=\tau(C) \text{ for all } k \geq q\}.$

\end{definition}

\begin{cor}
\label{cor2}
Let $C:f=0$ be a reduced plane curve of degree $d\geq 3$, let $r=mdr(f)$  be the minimal degree of a non-zero syzygy in  $AR(f)$. If $C$ is a maximal Tjurina curve, then the minimal resolution of the corresponding graded Milnor algebra $M(f)$, regarded as an $S$-module, has the form
$$0 \to S(-d-r)^{m-2} \to S(1-r-d)^m \to S(1-d)^3 \to S \to M(f) \to 0.$$
In particular, one has
$$ct(f)=st(f).$$

\end{cor}

\proof
The first claim follows from the obvious exact sequence
$$ 0 \to AR(f)(1-d) \to S(1-d)^3 \to S \to M(f) \to 0,$$
using Theorem \ref{thm2}. For the second claim, note that we have
$$\dim M(f)_k=\dim S_k-3 \dim S_{k+1-d}+m  \dim S_{k+1-d-r}-(m-2)  \dim S_{k-d-r}=$$
$$={k+2 \choose 2}-3{k+3-d \choose 2}+m {k+3-d-r \choose 2}-(m-2){k+2-d-r \choose 2},$$
if and only if $k \geq d+r-2$, where the binomial coefficients are regarded as polynomials in $k$ given by the usual formulas. It follows that $st(f)=d+r-2$.
On the other hand, it is known that 
$$ct(f)=d-2+mdr'(f),$$
 where $mdr'(f)$ is the minimal degree of a syzygy in $AR(f)$
which is not in the submodule $KR((f) \subset AR(f)$ generated by the Koszul relations $(f_y,-f_x,0)$, $(f_z,0,-f_x)$ and $(0,f_z,-f_y)$, see
\cite[Formula (1.3)]{DBull}. If $r<d-1$, then clearly $mdr'(f)=r$ and the last claim is proved.
If $r=d-1$ and $C$ is a maximal Tjurina curve, then $AR(f)$ is generated by
$2r-d+3=d+1 >3$ elements, so at least one of them is not in the 3-dimensional vector space
$KR(f)_{d-1}$. This implies again $mdr'(f)=r$ and the last claim is proved in this case also.
\endproof

\begin{rk}
\label{rkctst}  There are curves $C:f=0$ which are not maximal Tjurina curves, but which satisfy the equality $ct(f)=st(f).$ To see this, consider a uninodal curve $C$ of degree $d$, as in Remark \ref{rkS1} (iii) above, for which
$$ct(f)=st(f)=3(d-2).$$
When $d >3$, such a curve is not  a maximal Tjurina curve.
\end{rk}

\section{Existence of maximal Tjurina curves  when $2r \geq d$} 

Our discussion is a case-by-case analysis, according to the positive integer $$m=2r-d+3.$$
In the first two subsections we consider small values of $m$.
\subsection{Maximal Tjurina curves in the case $m=3$, minimal value for $m$} \label{s1}
This corresponds to the case $d=2r$ even, and the global Tjurina number is given by
$$\tau(C)=(d-1)(d-r-1)+r^2-1.$$
Then it follows from \cite{Dmax} that this equality occurs exactly when $C$ is a nearly free curve, with exponents $d_1=d_2=d_3=r$. Examples of such nearly free curves, both irreducible and line arrangements, are given in \cite{DStExpo}, for any pair $(d,r)=(2r,r)$.
From now we will assume $2r>d$.

\subsection{Maximal Tjurina curves in the case $m=4$}  \label{s2}
In this case the degree $d=2r-1$ is odd, and according to Theorem \ref{thm2}, the exponents are $d_1=d_2=d_3=d_4=r$. Such curves have occurred in \cite[Theorem 3.11]{DSt3syz}, 
and examples for the pairs $(d,r) \in \{ (5,3), (7,4), (9,5)\}$ are given in \cite[Example 3.12]{DSt3syz}.

The following example gives a sequence of maximal Tjurina curves which are in the same time {\it rational nearly cuspidal curves}.

\begin{prop}
\label{prop0}  Let $d =2r-1\geq 5$ be an odd integer and set 
$$C_d:f_d=(y^3 - x^2z)x^{r-3}y^{r-1}+x^d+y^d=0.$$
Then the plane curve $C_d$ is a maximal Tjurina curve of type $(d,r)$ for any odd degree $d\geq 5$. Moreover, any  curve $C_d$ is rational, has a unique singular point, namely $p=(0:0:1)$, and the plane curve singularity $(C_d,p)$ has two branches. 
\end{prop}

\proof The minimal degree syzygy for $f_d$ is given by
$$\rho_1=(0,x^{r-1}y, (r+2)x^{r-3}y^3 + (2r-1)y^r -
(r-1)x^{r-1}z),$$ 
and hence indeed $mdr(f_d)=r$. 
The curve $C_d$ is clearly rational, since we can express $z$ as a rational function of $x$ and $y$. The Milnor number $\mu(C_d,p)$ can be easily computed, since the singularity $(C_d,p)$ is  Newton nondegenerate and commode, see \cite{K}. It follows that
$$\mu(C_d,p)=4r^2-10r+5.$$
Since $C_d$ is rational, we have for the $\delta$-invariant the following equality
$$\delta(C_d,p)=\frac{(d-1)(d-2)}{2}=(r-1)(2r-3).$$
It follows that the number of branches of the singularity $(C_d,p)$ is 
$$2\delta(C_d,p) -\mu(C_d,p)+1=2,$$
and hence $C_d$ is a nearly cuspidal rational curve. Apply now \cite[Theorem 5.5]{DSt3syz}
with $d'=r-1$, and conclude that $C$ is a 4-syzygy curve with exponents
$d_1=d_2=d_3=d_4=r$. The last claim in Theorem \ref{thm2} implies that $C$ is indeed a maximal Tjurina curve of type $(d,r)$.
\endproof

Next we construct maximal Tjurina line arrangements for $m=4$.

\begin{ex}
\label{exm=4}  For $r \geq 4$ consider the line arrangement
$$\A_d: f=(x-z)(x-2z)…(x-(r-2)z)(y-z)(y-2z)…(y-(r-2)z)z(y-x-z)(y-x-2z)=0,$$
consisting of $d=2r-1$ lines. It is easy to check that 
$$\tau(\A_d)=3r^2-6r+1,$$
and using \cite[Theorem 1.2]{Dcurves}, it follows that 
$d_1=mdr(f) \in \{r-1,r\}$, since $\A_d$ has two points of multiplicity $r-1$. Hence the claim that $\A_d$ is a maximal Tjurina curve of type $(2r-1,r)$ is equivalent to showing $mdr(f)=r$. We can check the equality  $mdr(f)=r$ only using the SINGULAR software, for all $r$ with $4 \leq r \leq 15$. Takuro Abe has a proof of the equality $mdr(f)=r$ in full generality.
 For the case of an arrangement of 5 lines, see the arrangement $\A_5$ in \cite[Example 3.12]{DSt3syz}.
\end{ex}
Now we consider large values of $m$, in decreasing order.

\subsection{Maximal Tjurina curves in the case $m=d+1$, maximal value for $m$} 

In this case $r=d-1$. The following examples have been checked using SINGULAR. 

\begin{ex}
\label{ex1}  Let $d=2p \geq 4$ be an even  integer and set
$$f=(x^2-yz)^{p-1}yz+x^d+y^d.$$
Then the plane curve $C:f=0$ is a maximal Tjurina curve of type $(d,d-1)$ for any even degree $d=2p$ with $2 \leq p \leq 15$.
\end{ex}

\begin{ex}
\label{ex2}  Let $d =2p+1\geq 5$ be an odd integer and set
$$f=(x^2-yz)^{p-1}xyz+x^d+y^d.$$
Then  the plane curve $C:f=0$ is a maximal Tjurina curve of type $(d,d-1)$ for any odd degree $d=2p+1$ with $2 \leq p \leq 15$.
\end{ex}

Recall that for any $d \geq 2$ there are irreducible, rational, nodal curves of degree $d$. They  have exactly
$(d-1)(d-2)/2$ nodes and no other singularities, see \cite{Har, Oka}. For these curves, which are called {\it maximal nodal curves} in \cite{Oka}, we have the following result.

\begin{prop}
\label{prop1}  Let $C$ be a maximal nodal curve of degree $d$. Then $C$ is a maximal Tjurina curve of type $(d,d-1)$.
\end{prop}

\proof
First note that $\tau(d,d-1)_{max}=(d-1)(d-2)/2$, hence it remains to recall that an irreducible
nodal curve $C:f=0$ has $d_1=mdr(f)=d-1$, see \cite[Theorem 4.1]{DStEdin}.

\endproof

\begin{rk}
\label{rkLA}  If $\A:f=0$ is an arrangement of $d>1$ lines, with a point of maximal multiplicity $m(\A)$, then it follows from \cite[Theorem 1.2]{Dcurves} that $mdr(f) \leq d-m(\A) \leq d-2$. Hence there are no maximal Tjurina line arrangements of type $(d,d-1)$ when $d>1$.
\end{rk}

\subsection{Maximal Tjurina curves in the case $m=d-1$.} 
In this case $r=d-2$ and here are some examples.
$$(d,r)=(6,4)   \text{ and } f=(y^2z - x^3)^2 +x^6+y^6+xy^5.$$
$$(d,r)=(7,5)  \text{ and } f=(y^2z - x^3)^2y +x^7+y^7.$$
$$(d,r)=(8,6)  \text{ and } f=(y^2z - x^3)^2xy +x^8+y^8.$$
$$(d,r)=(9,7)  \text{ and } f=(y^3z+x^4)(x^3z+y^4)y+x^9+y^9.$$
$$(d,r)=(10,8)  \text{ and } f= (y^2z - x^3+x^2y)^3y +x^{10}+y^{10}.$$
The fact that these curves are maximal Tjurina curves can be checked using a computer algebra software, for instance SINGULAR.
In this case $r=d-2> d/2$ implies $d\geq 5$, and a direct computation shows that
$$\tau(d,d-2)_{max}={d \choose 2}.$$
One has the following result about maximal Tjurina line arrangements in $\PP^2$ of type $(d,d-2)$.
\begin{prop}
\label{prop0.1}  Let $C$ be a generic arrangement of $d\geq 4 $ lines in $\PP^2$. Then $C$ is a maximal Tjurina curve of type $(d,d-2)$.
\end{prop}

\proof
First note that $\tau(d,d-2)_{max}=\tau(C)$, since $C$ has only nodes as singularities and their number is given by ${d \choose 2}.$ It remains to recall that any reducible
nodal curve $C:f=0$ has $d_1=mdr(f)=d-2$, see \cite[Theorem 4.1]{DStEdin}.

\endproof

\subsection{Maximal Tjurina curves in the case $m=d-3$.} \label{s3}
In this case $r=d-3> d/2$ implies $d \geq 7$, and a direct computation shows that
$$\tau(d,d-3)_{max}={d+1 \choose 2}-3.$$
We construct a sequence of line arrangements $C_d:f_d=0$ such that
$C_d$ consists of $d$ lines, has only double and triple points, 
$\tau(C_d)=\tau(d,d-3)_{max}$ and $r=mdr(f_d)=d-3$.
We consider two sequences $a_n=2^n$ and $b_n=3^n$. For each integer $e\geq 2$ we set
$$h_{2e-2}=\prod_{j=1}^{j=e-1}\left(\frac{x}{a_j}+\frac{y}{b_j}-z\right)\left(\frac{x}{a_j}+\frac{y}{b_{j+1}}-z\right).$$
Then we define our polynomials
$$f_{2e}=xyh_{2e-2} \text{ and } f_{2e+1}=xy\left(\frac{x}{a_e}+\frac{y}{b_e}-z\right)h_{2e-2}.$$
Then we conjecture that the arrangement $C_d:f_d=0$ is Tjurina maximal of type $(d,d-3)$ for $d \geq 7$. The claim $\tau(C_d)=\tau(d,d-3)_{max}$  is easy to check, just by counting the number of double and triple points.
The claim $mdr(f_d)=d-3$ is more difficult to check. Indeed, Ziegler's celebrated example of two arrangements $\A:f=0$ and $\A': f'=0$, both consisting of 9 lines and having only double and triple points, with isomorphic intersection lattices, and $mdr(f)=6$ and $mdr(f')=5$, shows that the invariant $mdr(f)$ is not combinatorial, see \cite{Zi} and \cite[Remark 8.5]{DHA}. We have checked the claim
 $mdr(f_d)=d-3$ for all degrees $d$ with $7 \leq d \leq 20$, using SINGULAR.

\subsection{Maximal Tjurina curves in the case $m=d-5$.} \label{s4}
In this case $r=d-4 > d/2$ implies $d \geq 9$, and a direct computation shows that
$$\tau(d,d-4)_{max}={d+2 \choose 2}-9.$$
We construct a sequence of line arrangements $\A_d:g_d=0$ such that
$\A_d$ consists of $d$ lines, has only double, triple and 4-fold points, 
$\tau(\A_d)=\tau(d,d-4)_{max}$ and, conjecturally, $r=mdr(g_d)=d-4$.
To do this we use the sequences $a_n$, $b_n$ and the polynomials $f_{2e+1}$ constructed in the previous section, and define new polynomials as follows, for $k\geq 2$.
$$g_{3k+2}=f_{2k+3} \prod_{p=1}^{p=k-1}\left(\frac{x}{a_p}+\frac{y}{b_{p+2}}-z\right),$$
$$g_{3k+3}=(27x-8y)g_{3k+2} \text{ and }g_{3k+4}=(x-y)g_{3k+3}.$$
The claim $\tau(\A_d)=\tau(d,d-4)_{max}$  is easy to check, just by counting the number of double, triple and 4-fold points. We have checked the claim
 $mdr(g_d)=d-4$ for all degrees $d$ with $9 \leq d \leq 20$, using SINGULAR.

\end{document}